\documentclass{gtmon_a}    
\pdfoutput=1
\usepackage{amscd}

\proceedingstitle{Groups, homotopy and configuration spaces (Tokyo
  2005)}
\conferencestart{5 July 2005}
\conferenceend{11 July 2005}
\conferencename{Groups, homotopy and configuration spaces,
                in honour of Fred Cohen's 60th birthday}
\conferencelocation{University of Tokyo, Japan}

\editor{Norio Iwase}
\givenname{Norio}
\surname{Iwase}

\editor{Toshitake Kohno}
\givenname{Toshitake}
\surname{Kohno}

\editor{Ran Levi}
\givenname{Ran}
\surname{Levi}

\editor{Dai Tamaki}
\givenname{Dai}
\surname{Tamaki}

\editor{Jie Wu}
\givenname{Jie}
\surname{Wu}

\title{Twisted Morita--Mumford classes on braid groups}

\author{Nariya Kawazumi}
\givenname{Nariya}
\surname{Kawazumi}
\address{
Department of Mathematical Sciences\\
University of Tokyo\\\newline
Tokyo, 153-8914 Japan 
}                  
\email{kawazumi@ms.u-tokyo.ac.jp}                     
\urladdr{}

\keyword{free group}
\keyword{Morita-Mumford class}
\keyword{automorphism group of a free group}
\keyword{braid group}
\keyword{cohomology of groups}
\keyword{mapping class group}

\subject{primary}{msc2000}{20F36}
\subject{secondary}{msc2000}{14H15}
\subject{secondary}{msc2000}{20J06}
\subject{secondary}{msc2000}{20F28}
\subject{secondary}{msc2000}{32G15}
\subject{secondary}{msc2000}{57R20}
\subject{secondary}{msc2000}{57M50}

\arxivreference{math.GT/0606102}

\volumenumber{13}
\issuenumber{}
\publicationyear{2008}
\papernumber{13}
\startpage{293}
\endpage{306}
\doi{}
\MR{}
\Zbl{}

\received{5 June 2006}
\revised{28 May 2007}
\accepted{30 May 2007}
\published{25 February 2008}
\publishedonline{25 February 2008}
\proposed{}
\seconded{}
\corresponding{}
\version{}

\arxivreference{}

\makeatletter
\def\cnewtheorem#1[#2]#3{\newtheorem{#1}{#3}[section]
\expandafter\let\csname c@#1\endcsname\c@thm}

\AtBeginDocument{\let\bar\wbar\let\tilde\wtilde\let\hat\what\let\Bbb\mathbb}

\newtheorem{theorem}{Theorem}
\newtheorem{thm}{Theorem}[section]  
\cnewtheorem{lem}[thm]{Lemma}        

\cnewtheorem{prop}[thm]{Proposition}   %
\cnewtheorem{cor}[thm]{Corollary}   %
\cnewtheorem{expect}[thm]{Expectation}   %
\theoremstyle{definition}
\cnewtheorem{defn}[thm]{Definition}    

\makeatother  

\numberwithin{equation}{section}

\newcommand\abel{^{\mbox{abel}}}
\newcommand\An{{\mbox{Aut}(F_n)}}
\newcommand\Anbar{{\overline{A_n}}}
\newcommand\bb{{\varsigma}}
\newcommand\bC{{\bf C}}
\newcommand\bZ{{\bf Z}}
\newcommand\bQ{{\bf Q}}

\newcommand\bbn{{\widehat{B_n}}}
\newcommand\bn{{B_n}}
\newcommand\Fn{{F_n}}
\newcommand\h{{\overline{h}}}
\newcommand\Hom{{\mbox{Hom}}}
\newcommand\HQ{{H_{\bf Q}}}
\newcommand\HcQ{{H^*_{\bf Q}}}
\newcommand\IAn{{IA_n}}
\newcommand\inv{^{-1}}

\newcommand\La{{\bigwedge}}
\newcommand\MCG{{\mathcal{M}_{g, 1}}}
\newcommand\pnq{{\mathcal{P}_{n-q}(q)}}
\newcommand\sn{{\mathfrak{S}_n}}
\newcommand\std{{\mbox{std}}}
\newcommand\T{{\widehat{T}}}
\newcommand\xxi{{\widehat{\xi}}}

\begin{document}

\begin{htmlabstract}
Evaluating the twisted Morita&ndash;Mumford classes <span
style="text-decoration: overline">h</span><sub>p</sub> on
the Artin braid group B<sub>n</sub>, we give the stable algebraic independence
of the <span style="text-decoration: overline">h</span><sub>p</sub>'s on the automorphism group of the free group,
Aut(F<sub>n</sub>).  This is sharper than the results we obtained by
restricting them to the mapping class group.
\end{htmlabstract}

\begin{abstract} 
Evaluating the twisted Morita--Mumford classes $\overline{h}_p$
(Kawazumi \cite{Ka4}) on the Artin braid group $B_n$, we give the
stable algebraic independence of the $\overline{h}_p$'s on the
automorphism group of the free group, $\mbox{Aut}(F_n)$.  This is
sharper than the results obtained by restricting them to the mapping
class group (Kawazumi \cite{Ka2}).
\end{abstract}
\begin{webabstract} 
Evaluating the twisted Morita--Mumford classes $\overline{h}_p$ on
the Artin braid group $B_n$, we give the stable algebraic independence
of the $\overline{h}_p$'s on the automorphism group of the free group,
$\mathrm{Aut}(F_n)$.  This is sharper than the results we obtained by
restricting them to the mapping class group.
\end{webabstract}
\begin{asciiabstract} 
Evaluating the twisted Morita-Mumford classes bar(h)_p on the Artin
braid group B_n, we give the stable algebraic independence of the
bar(h)_p's on the automorphism group of the free group, Aut(F_n).
This is sharper than the results we obtained by restricting them to
the mapping class group.
\end{asciiabstract}

\maketitle


\section*{Introduction}

In the cohomological study of the mapping class group for a surface,
the Morita--Mumford classes, $e_i = (-1)^{i+1}\kappa_i$, $i \geq 1$,
\cite{Mu,Mo1} play some important roles. As was proved by Miller
\cite{Mi} and Morita \cite{Mo1} independently, they are algebraically
independent in the stable range $* < \frac23 g$.  Madsen and Weiss
\cite{MW} proved that the rational stable cohomology algebra of the
mapping class groups, $H^*(\mathcal{M}_\infty; {\bf Q})$, is generated
by the Morita--Mumford classes.  The Morita--Mumford classes have
twisted variants, $m_{i, j} \in H^{2i+j-2}(\MCG; \La^jH)$, $i, j \geq
0$, introduced by the author \cite{Ka1}.  Here we denote by
$\Sigma_{g, 1}$ a $2$--dimensional oriented compact connected
$C^\infty$ manifold of genus $g$ with $1$ boundary component, $\MCG$
its mapping class group, $\MCG := \pi_0\mbox{Diff}(\Sigma_{g, 1},
\mbox{id on $\partial \Sigma_{g, 1}$})$, and $H$ the integral first
homology group of the surface $\Sigma_{g, 1}$.  The mapping class
group $\MCG$ acts on $H$ in an obvious way.  The twisted\break variants also
satisfy the algebraic independence.  More precisely, the algebra\break
$H^*(\MCG; \La^*H)\otimes \bQ$ is the polynomial algebra in the set
$\{m_{i, j}; \, i \geq 0, j \geq 1, \,\mbox{and}\, i+j \geq 2\}$ over
the algebra $H^*(\MCG; \bQ)$ in the range where the total degree $\leq
\frac23 g$ (Kawazumi \cite[Theorem 1.C]{Ka2}.) Hence, from the theorem
of Madsen and Weiss \cite{MW} stated above, the algebra $H^*(\MCG;
\La^*H)\otimes \bQ$ is stably isomorphic to the polynomial algebra in
the set $\{m_{i, j}; \, i \geq 0, j \geq 0, \,\mbox{and}\, i+j \geq
2\}$ over $\bQ$.  Similar results hold for any other symplectic
coefficients (Kawazumi \cite[Theorem 1.B]{Ka2}.)  Furthermore all the cohomology
classes on the mapping class group obtained by contracting the
coefficients of the twisted ones using the intersection pairing
$H^{\otimes 2} \to \bZ$ are exactly the algebra generated by the
(original) Morita--Mumford classes $e_i$'s (Morita \cite{Mo2}, Kawazumi
and Morita \cite{KM}).

Some of the twisted ones have the advantage over the original ones of
being defined on the automorphism group of a free group, which has the
mapping class group and the braid group as proper subgroups.  Let $n
\geq 2$ be an integer, $\Fn$ a free group of rank $n$ with free basis
$x_1, x_2, \dots, x_n$
$$
\Fn = \langle x_1, x_2, \dots, x_n\rangle,
$$
and $\An$ the automorphism group of the group $\Fn$. 
The Dehn--Nielsen theorem tells us the natural action of 
the group $\MCG$ on the free group $\pi_1(\Sigma_{g, 1})$ 
of rank $2g$ induces an injective homomorphism $\MCG \to 
\mbox{Aut}(F_{2g})$. In view of a theorem of Artin 
\cite{A} the braid group $\bn$ of $n$ strings is embedded into 
the group $\An$. \par
Now we denote by $H$ and $H^*$ the first integral homology 
and cohomology groups of the group $\Fn$
$$
H := H_1(\Fn; \bZ) = \Fn\abel = \Fn/[\Fn. \Fn] \quad\mbox{and}\quad
H^* := H^1(\Fn; \bZ) = \Hom(H, \bZ), 
$$ 
respectively, on which the automorphism group $\An$ acts in 
an obvious way. We write $[\gamma] := \gamma \bmod [\Fn, \Fn] \in H$
for $\gamma \in \Fn$, and $X_i := [x_i] \in H$  for $i$, $1 \leq
i\leq n$. 
In \cite{Ka4} we introduced cohomology classes 
$$
h_p \in H^p(\An; H^*\otimes H^{\otimes (p+1)})\quad
\mbox{and}\quad \h_p \in H^p(\An; H^{\otimes p})
$$
for $p \geq 1$. Restricted to the mapping class group 
$\MCG$ they coincide with the twisted Morita--Mumford classes
\begin{eqnarray*}
&& (p+2)!\,h_p\vert_{\MCG} = m_{0, p+2} \in H^p(\MCG; H^{\otimes (p+2)}),
\quad\mbox{and}\\
&& p!\,\h_p\vert_{\MCG} = -m_{1, p} \in H^p(\MCG; H^{\otimes p}).
\end{eqnarray*}
Here $H$ and $H^*$ are isomorphic to each other as $\MCG$ modules 
because of the intersection pairing of the surface $\Sigma_{g,1}$. 
The class $p!\h_p$ can be regarded as an element in 
$H^p(\An; \La^pH)$. \par
In this note we confine ourselves to studying the behavior of 
$\h_p$'s restricted to the braid group $\bn$, and consider 
the rational coefficients
$$
\HQ := H\otimes_{\bZ}\bQ \quad \mbox{and}\quad
\HcQ := H^*\otimes_\bZ\bQ.
$$ 
In this paper we prove the following result:

\begin{theorem}\label{0.1}
The cohomology classes $\h_p$'s are algebraically independent 
in the algebra $H^*(\bn; \La^*\HQ)$ in the range where 
the total degree $\leq n$.
\end{theorem}

Here the total degree of $\h_p$ is defined to be $2p$. 
\fullref{0.1} implies the algebraic independence on the automorphism 
group $\An$. This is sharper than that obtained by restricting 
them to the mapping class group $\MCG$ \cite[Theorem 1.C]{Ka2}, 
where the range is given by the inequality the total degree 
$\leq \frac23 g = \frac13 n$. 
\par
\fullref{0.1} was announced in \cite{Ka3}. Its proof given in \fullref{sec3} is based on
some kind of primitiveness of the $\h_p$'s (\fullref{1.2}) and the
evaluation of $\h_{n-1}$ on the pure braid group of $n$ strings, $P_n$
(\fullref{2.4}).
In \fullref{sec4} we will give some remarks on the cohomology of 
the automorphism group $\An$.

\section{Twisted Morita--Mumford classes on the automorphism\break
group $\An$}\label{sec1}

Throughtout this paper we denote by $C^*(G; M)$ 
the normalized standard complex of a group $G$ with values in 
a $G$--module $M$, and use the Alexander--Whitney cup product 
$\cup\co  C^*(G; M_1) \otimes C^*(G; M_2)\to C^*(G; M_1\otimes
M_2)$. Moreover we denote by $Z^p(G; M)$, $p \geq 0$, the
$p$--cocycles in  the cochain complex $C^*(G; M)$. \par

Now we recall the definition of the twisted cohomology classes $h_p$ 
and $\h_p$ on the automorphism group $\An$ for $p \geq 1$. 
The semi-direct product 
$$
\Anbar := \Fn\rtimes\An
$$ 
admits an extension of groups
\begin{equation}\label{1-1}
\Fn \mathop\to^\iota \Anbar \mathop\to^\pi \An
\end{equation}
given by $\iota(\gamma) = (\gamma, 1)$ and $\pi(\gamma, \varphi) =
\varphi$  for $\gamma \in \Fn$ and $\varphi \in \An$. 
The map $k_0: \Anbar \to H$, $(\gamma, \varphi) \mapsto [\gamma]$,
satisfies the cocycle condition. 
We write also $k_0$ for the cohomology class 
$[k_0] \in H^1(\Anbar; H)$. 
For each $p \geq 1$ we define $h_p$ by the image of the 
$(p+1)$-st power of the cohomology class $k_0$ 
under the Gysin map of the extension \eqref{1-1}
\begin{equation}
h_p := \pi_\sharp({k_0}^{\otimes (p+1)}) \in H^p(\An; 
H^*\otimes H^{\otimes (p+1)})
\end{equation}\label{1-2}
\cite{Ka4}. Contracting the coefficients by the
$\mbox{GL}(H)$--homomorphism
\begin{equation}\label{1-3}
r_p\co H^*\otimes H^{\otimes (p+1)} 
\to H^{\otimes p},\quad
f\otimes v_0\otimes v_1 \otimes \cdots \otimes v_p \mapsto 
f(v_0) v_1 \otimes \cdots \otimes v_p,
\end{equation}
we define
\begin{equation}\label{1-4}
\h_p := {r_p}_*(h_p) \in H^p(\An; H^{\otimes p}).
\end{equation}
The $p$-th exterior power ${k_0}^p = p!{k_0}^{\otimes p}$ 
can be regarded as a cohomology class with coefficients in 
$\La^pH$. Hence, if we consider the rational coefficients 
$\HQ$, we may regard $\h_p$ as a cohomology class in
$H^p(\An; \La^p\HQ)$.  
\par
A Magnus expansion $\theta$ of the free group $\Fn$ 
gives an explicit cocycle representing the class $h_p$. 
The completed tensor algebra generated by $H$, 
$\T = \T(H) : = {\prod}^\infty_{m=0} H^{\otimes m}$, 
has a decreasing filtration of two-sided ideals
$\T_p := {\prod}_{m\geq p} H^{\otimes m}$, $p \geq 1$.
It should be remarked that the subset $1+\T_1$ is 
a subgroup of the multiplicative group of the algebra $\T$.
We call a map $\theta\co  \Fn \to 1 + \T_1$ a {\it Magnus expansion} 
of the free group $\Fn$, if $\theta\co  \Fn \to 1 + \T_1$ is a group
homomorphism, and if $\theta(\gamma) \equiv 1 + [\gamma] \pmod{\T_2}$ 
for any $\gamma \in \Fn$. 
We write $\theta(\gamma) = \sum^\infty_{m=0}\theta_m(\gamma)$, 
$\theta_m(\gamma) \in H^{\otimes m}$. The $m$-th component 
$\theta_m\co  \Fn \to H^{\otimes m}$ is a map, but {\it not} 
a group homomorphism. 
A Magnus expansion $\std\co  \Fn \to 1+\T_1$ is 
defined by $\std(x_i) := 1 + X_i$, $1 \leq i \leq n$.
Here we denote $X_i := [x_i] \in H$, the homology class 
of the generator $x_i$. 
We call it {\it the standard Magnus expansion}. 
As is described in classical references, 
the value $\std(\gamma)$ for any word $\gamma \in \Fn$ 
is explicitly computed by means of Fox' free differentials. 
All the results of this paper can be derived from 
the expansion $\std$. \par

We define a map 
$\tau^\theta_1\co  \An \to H^*\otimes H^{\otimes 2}$
by 
\begin{equation}\label{1-5}
\tau^\theta_1(\varphi)[\gamma] 
= \theta_2(\gamma) 
-\vert\varphi\vert^{\otimes 2}\theta_2(\varphi\inv(\gamma)) \in
H^{\otimes 2}
\end{equation}
for $\gamma \in \Fn$ and $\varphi \in \An$. 
Here $\vert\varphi\vert \in \mbox{GL}(H)$ is the 
automorphism of $H = \Fn\abel$ induced by $\varphi$. 
This map $\tau^\theta_1$ satisfies the cocycle condition 
\cite[Lemma 2.1]{Ka4}. Now we introduce a $\mbox{GL}(H)$--homomorphism
$$
\bb_p\co  (H^*\otimes H^{\otimes 2})^{\otimes p} = 
\Hom(H, H^{\otimes 2})^{\otimes p} \to \Hom(H, H^{\otimes (p+1)})
= H^*\otimes H^{\otimes (p+1)}
$$
for each $p \geq 1$. If $p \geq 2$, we define
\begin{eqnarray}\label{1-6}
&& \bb_p(u_{(1)}\otimes u_{(2)}\otimes
\cdots\otimes u_{(p-1)}\otimes u_{(p)}) \\
& := &
\left(u_{(1)}\otimes {1_H}^{\otimes (p-1)}\right)\circ
\left(u_{(2)}\otimes {1_H}^{\otimes (p-2)}\right)\circ\cdots\circ
\left(u_{(p-1)}\otimes {1_H}\right)\circ u_{(p)},\nonumber
\end{eqnarray}
where $u_{(i)} \in \Hom(H, H^{\otimes 2}) = H^*\otimes H^{\otimes 2}$, 
$1 \leq i \leq p$.  
In the case $p=1$, we define $\bb_1 := 1_{H^*\otimes H^{\otimes 2}}$.
Then we have:
\begin{thm}{\rm \cite[Theorem 4.1]{Ka4}}\label{1.1} 
$$
h_p = {\bb_p}_*([\tau^\theta_1]^{\otimes p}) 
\in H^p(\An; H^*\otimes H^{\otimes(p+1)})
$$
for any Magnus expansion $\theta$ and each $p \geq 1$. 
In the case $p=1$ we have $[\tau^\theta_1] = h_1
\in  H^1(\An; H^*\otimes H^{\otimes 2})$.
\end{thm}

Some kind of primitiveness of the cohomology classes $h_p$ and 
$\h_p$ follows from the theorem.
We write simply $A_n := \An$ for the remainder of the section. 
Suppose $n_1 + n_2 \leq n$. 
Let $A_{n_2}$ act on the words in the letters $x_{n_1+1}, x_{n_1+2}, 
$\linebreak $\dots, x_{n_1+n_2}$ in an obvious way. 
Then we have a natural homomorphism
$$
\iota = 
\iota_{n_1, n_2}\co  A_{n_1}\times A_{n_2} \to A_n.
$$
We denote by $\varpi_1\co  A_{n_1}\times A_{n_2} \to 
A_{n_1}$ and $\varpi_2\co  A_{n_1}\times A_{n_2} \to 
A_{n_2}$ the first and the second projections of the product 
$A_{n_1}\times A_{n_2}$, respectively, and 
by $H_{(n_1)}$, $H_{(n_2)}$ and $H_{(n-n_1-n_2)}$ 
the submodules of $H$ spanned by 
$\{X_{1}, \dots, X_{n_1}\}$, $\{X_{n_1+1}, \dots, X_{n_1+n_2}\}$ 
and $\{X_{n_1+n_2+1}, \dots, X_{n}\}$, respectively. 
Then we have a direct-sum decomposition $H = H_{(n_1)}\oplus
H_{(n_2)}\oplus H_{(n-n_1-n_2)}$, and can consider the map
$$
{\varpi_k}^*\co  H^*(A_{n_k}; H^*_{(n_k)}\otimes H^{\otimes (p+1)}_{(n_k)}) 
\to H^*(A_{n_1}\times A_{n_2}; H^*\otimes H^{\otimes (p+1)})
$$
for $k = 1$ and $2$. For any $p \geq 1$ we have:

\begin{prop}\label{1.2}\

\begin{enumerate}
\item $\iota^*h_p = {\varpi_1}^*h_p + {\varpi_2}^*h_p \in
H^p(A_{n_1}\times A_{n_2}; H^*\otimes H^{\otimes (p+1)})$,
\item $\iota^*\h_p = {\varpi_1}^*\h_p + {\varpi_2}^*\h_p \in
H^p(A_{n_1}\times A_{n_2}; H^{\otimes p})$.
\end{enumerate}
\end{prop}

\begin{proof}
Using the standard expansion $\std$, we write simply 
$$\tau^{(k)} := {\varpi_k}^*\tau^\std_1 
\in Z^1(A_{n_1}\times
A_{n_2}; H^*\otimes H^{\otimes 2}).$$ 
Clearly we have $\std(\gamma_1) \in 
\prod^\infty_{p=0} {H_{(n_1)}}^{\otimes p} 
\subset \T$ for any word $\gamma_1$ in the letters 
${x_{1}}, \dots, {x_{n_1}}$. 
Similar conditions hold for any word $\gamma_2$
in the letters ${x_{n_1+1}}, \dots, {x_{n_1+n_2}}$
and any $\gamma_3$ in 
${x_{n_1+n_2+1}}, \dots, {x_{n}}$. 
Hence, from the definition of 
$\tau^\theta_1$ \eqref{1-5}, we have 
$$
\iota^*\tau^\std_1 = \tau^{(1)} + \tau^{(2)} \in 
Z^1(A_{n_1}\times
A_{n_2}; H^*\otimes H^{\otimes 2}).
$$
If we use the $\mbox{GL}(H)$--homomorphism 
${\bb_2}\co  (H^*\otimes H^{\otimes 2})^{\otimes 2} \to 
H^*\otimes H^{\otimes 3}$ in \eqref{1-6}, then 
we have
\begin{equation}\label{1-7}
{\bb_2}_*(\tau^{(1)}\tau^{(2)}) = {\bb_2}_*(\tau^{(2)}\tau^{(1)}) = 0 \in
Z^2(A_{n_1}\times A_{n_2}; H^*\otimes H^{\otimes 3}).
\end{equation}
In fact, $f(u) = 0$ for any $f \in H^*_{(n_1)}$ and 
$u \in H_{(n_2)}$ and vice versa. From \fullref{1.1} follows
\begin{eqnarray*}
& \iota^*h_p = {\bb_p}_*(\iota^*[\tau^\std_1]^{\otimes p}) 
= {\bb_p}_*((\tau^{(1)} + \tau^{(2)})^{\otimes p}) \\
& = {\bb_p}_*((\tau^{(1)})^{\otimes p}) + {\bb_p}_*((\tau^{(2)})^{\otimes
p})  = {\varpi_1}^*h_p + {\varpi_2}^*h_p.
\end{eqnarray*}
Here ${\bb_p}_*$ of each mixed term in $\tau^{(1)}$ and
$\tau^{(2)}$ vanishes by \eqref{1-7}. Applying ${r_p}_*$ to (1), 
we deduce (2). This completes the proof of the proposition.
\end{proof}

\section{Evaluation on the Artin braid groups}\label{sec2}

The $n$-th symmetric group 
$\sn$ acts on the space $\bC^n$ by permuting the components. 
The open subset 
$$
Y_n := \{(z_1, z_2, \dots, z_n) \in \bC^n;\,\,
 z_i \neq z_j \,\,\mbox{for $i \neq j$}\}
$$
is stable under the action of the group $\sn$. By definition, the  
Artin braid group of $n$ strings, $\bn$, is the fundamental group 
of the quotient space 
$Y_n/\sn$, $\bn := \pi_1(Y_n/\sn)$. As was shown by Artin \cite{A}, 
the group $B_n$ admits a presentation
\begin{eqnarray}\label{2-1}
\mbox{generators: }& \sigma_i, \quad 1 \leq i \leq n-1,\nonumber\\
\mbox{relations: }& \sigma_i\sigma_j = \sigma_j\sigma_i, \quad \mbox{if $\vert i
  - j\vert \geq 2,$}\\
& \sigma_i\sigma_{i+1}\sigma_i = \sigma_{i+1}\sigma_i\sigma_{i+1}, \quad
\mbox{for $1 \leq i \leq n-2$.}\nonumber
\end{eqnarray}
The pure braid group of $n$ strings, $P_n$, is defined to
be the fundamental group of the space $Y_n$, $P_n := \pi_1(Y_n)$. 
We have a natural extension of groups
$$
P_n \to B_n \to \sn. 
$$ 
As is known, $A_{i, j}$, $1 \leq i < j \leq n$, given by
$$
A_{i, j} := {\sigma_{j-1}}{\sigma_{j-2}}\cdots{\sigma_{i+1}}{\sigma_{i}}^2
{\sigma_{i+1}}\inv\cdots{\sigma_{j-2}}\inv{\sigma_{j-1}}\inv
$$
can serve as a generating system of the group $P_n$.
For details, see Birman \cite{B}. 
\par
The braid group $B_n$ admits a natural homomorphism 
into the group $\An$, $\xi\co  \bn \to \An$. 
To recall how to construct it, we consider an action of the 
group $\sn$ on the space $Y_{n+1} \subset 
\bC^{n+1} = \bC^n\times\bC$ given by 
$$\rho(z_1, \dots, z_n, z_{n+1}) = (z_{\rho^{-1}(1)}, 
\dots, z_{\rho^{-1}(n)}, z_{n+1})$$
for $\rho \in \sn$. 
We denote by $\bbn$ the fundamental group of the quotient space 
$Y_{n+1}/\sn$, 
$\bbn := \pi_1(Y_{n+1}/\sn)$. 

The forgetful map $Y_{n+1} \to Y_n$, $(z_1, \dots, z_n, z_{n+1}) 
\mapsto (z_1, \dots, z_n)$, induces a fibration
$$
\bC \setminus \{\mbox{$n$ points}\} \to Y_{n+1}/\sn \to Y_{n}/\sn
$$
with a section $s\co  Y_n/\sn \to Y_{n+1}/\sn$ given by
$(z_1, \dots, z_n) \mapsto (z_1, \dots, z_n,$ \linebreak
$\frac1n{\sum}^n_{i=1}z_i + {\sum}^n_{j=1}\vert z_j 
- \frac1n{\sum}^n_{i=1}z_i\vert)$
(Arnol'd \cite{Ar}). This fibration with the section $s$ induces 
an extension of groups
\begin{equation}\label{2-2}
\Fn \mathop\to^\iota \bbn \mathop\to^\pi \bn
\end{equation}
with a split homomorphism $s\co  \bn \to \bbn$. 
Thus we obtain a morphism of extensions of groups
\begin{equation}\label{2-3}
\CD
\Fn  @>>> \bbn @>>> \bn \\
@| @V{\xxi}VV @V{\xi}VV\\  
\Fn @>>> \Anbar @>>> \An.
\endCD
\end{equation}
The homomorphisms $\xi$ and $\xxi$ are explicitly given by
\begin{eqnarray*}
\iota(\xi(x)(\gamma)) &=& s(x)\gamma s(x)\inv\\
\xxi(\iota(\gamma)s(x)) &=& (\gamma, \xi(x)) \in \Fn\rtimes\An = \Anbar
\end{eqnarray*}
for $x \in \bn$ and $\gamma \in \Fn$. 
The group $\bbn$ is embedded into $B_{n+1}$ in an obvious way. 
Then the homomorphisms  
$s$ and $\iota$ are described as
\begin{eqnarray}\label{2-4}
s(\sigma_i) &=& \sigma_i \quad \mbox{for $1 \leq i \leq n -1$,}\\
\iota(x_j) &=& {\sigma_{n}}{\sigma_{n-1}}\cdots{\sigma_{j+1}}
{\sigma_{j}}^2{\sigma_{j+1}}\inv\cdots{\sigma_{n-1}}\inv{\sigma_{n}}\inv
\nonumber\\ 
 &=& A_{j, n+1}\quad \mbox{for $1 \leq j \leq n$}\nonumber
\end{eqnarray}
in terms of the presentation \eqref{2-1}. So the homomorphism $\xi$ is
explicitly given by
\begin{equation}\label{2-5}
\xi(\sigma_i)(x_j) = 
\begin{cases}
x_{i+1}, \quad& \mbox{if $j = i$,}\\
{x_{i+1}}\inv x_i x_{i+1}, \quad& \mbox{if $j = i+1$,}\\
x_j, \quad& \mbox{otherwise.}
\end{cases}
\end{equation}
We now evaluate the cohomology classes
$h_1$ and $\h_{n-1}$ on the braid group $\bn$. 
Here we use the standard Magnus expansion
$\std\co \Fn \to 1+\T_1$ introduced in \fullref{sec1}. 
For the rest of this section we write simply $k_0$, $\tau_1$, $h_p$ and
$\h_p$  for $\xxi^*k_0$, $\xi^*\tau^\std_1$, $\xi^*h_p$ and
$\xi^*\h_p$, respectively. 
Let $\{l_i\}^n_{i=1} \subset H^*$ denote the dual basis of 
 $\{X_i\}^n_{i=1} = \{[x_i]\}^n_{i=1} \subset H$. 

\begin{lem}\label{2.1}
$$
\tau_1(\sigma_i) 
= l_i\otimes (X_i \otimes X_{i+1} - X_{i+1} \otimes X_{i}) 
\in H^*\otimes H^{\otimes 2}
$$
\end{lem}
\begin{proof} From \eqref{1-5}
\begin{eqnarray*}
\tau_1(\sigma_i) 
&=& {\sum}^n_{j=1}l_j\otimes(\std_2(x_j) - \vert\sigma_i\vert^{\otimes 2}
\std_2({\sigma_i}\inv(x_j)))\\
&=& - l_i\otimes\vert\sigma_i\vert^{\otimes 2} \std_2({\sigma_i}\inv(x_i))
 - l_{i+1}\otimes\vert\sigma_i\vert^{\otimes 2}
\std_2({\sigma_i}\inv(x_{i+1}))\\
&=& - l_i\otimes\vert\sigma_i\vert^{\otimes 2} \std_2(x_ix_{i+1}{x_i}\inv)
 - l_{i+1}\otimes\vert\sigma_i\vert^{\otimes 2}
\std_2(x_i)\\
&=& - l_i\otimes\vert\sigma_i\vert^{\otimes 2}
\std_2(x_ix_{i+1}{x_i}\inv).
\end{eqnarray*}
On the other hand, we have 
$$
\std_2(x_ix_{i+1}{x_i}\inv) = X_i\otimes X_{i+1} - X_{i+1}\otimes X_i.
$$
In fact, $X_i\otimes X_{i+1} = \std_2(x_ix_{i+1}) 
= \std_2(x_ix_{i+1}{x_i}\inv x_i) 
= \std_2(x_ix_{i+1}{x_i}\inv) + \std_2(x_i) + X_{i+1}\otimes X_i
= \std_2(x_ix_{i+1}{x_i}\inv) + X_{i+1}\otimes X_i$. 
Therefore we obtain
$\tau_1(\sigma_i) =  - l_i\otimes\vert\sigma_i\vert^{\otimes
2}(X_i\otimes X_{i+1} - X_{i+1}\otimes X_i) = - l_i\otimes
(X_{i+1}\otimes X_i - X_i\otimes X_{i+1})$, as was to be shown.
\end{proof}

The pure braid group $P_n$ acts on the homology $H$ trivially. 
Hence, from \cite[Theorem 3.1]{Ka4}, the restriction of $\tau_1$ 
to $P_n$ does not depend on the choice of Magnus expansions. 

\begin{lem}\label{2.2}
$$
\tau_1(A_{i, j}) 
= (l_i - l_j)\otimes(X_i\otimes X_{j} - X_{j}\otimes X_{i})
$$
\end{lem}
\begin{proof} 
Recall the map $\tau_1$ satisfies the cocycle condition on the
automorphism group $\An$.  When we set $\gamma 
:= {\sigma_{j-1}}{\sigma_{j-2}}\cdots{\sigma_{i+1}}$, we have $A_{i, j} 
= \gamma{\sigma_{i}}^2\gamma\inv$, so that
\begin{eqnarray*}
&&\tau_1(A_{i, j}) \\
&=& \tau_1(\gamma{\sigma_{i}}^2\gamma\inv) 
= \tau_1(\gamma) + \gamma\tau_1({\sigma_{i}}^2) 
+ \gamma{\sigma_{i}}^2\tau_1(\gamma\inv)\\
&=& \tau_1(\gamma) + \gamma\tau_1({\sigma_{i}}^2) 
+ \gamma\tau_1(\gamma\inv) = \tau_1(1) + \gamma\tau_1({\sigma_{i}}^2) 
= \gamma\tau_1({\sigma_{i}}^2)\\ 
&=& \gamma(\tau_1({\sigma_{i}}) + {\sigma_{i}}\tau_1({\sigma_{i}}))\\
&=& \gamma(l_{i}\otimes(X_i\otimes X_{i+1} - X_{i+1}\otimes X_{i})) 
+ {\gamma\sigma_{i}}(l_{i}\otimes(X_i\otimes X_{i+1} - X_{i+1}\otimes
X_{i}))\\ &=& \gamma((l_{i} - l_{i+1})
\otimes(X_i\otimes X_{i+1} - X_{i+1}\otimes X_{i})) \\
&=& (l_{i} - l_{j})
\otimes(X_i\otimes X_{j} - X_{j}\otimes X_{i}),
\end{eqnarray*}
as was to be shown.
\end{proof}

To prove the nontriviality of $\h_{n-1}$ on the group $B_n$, 
we recall some basic facts on the cohomology of the pure braid group 
$P_n$. The space $Y_n$ is an Eilenberg--MacLane space of type 
$(P_n, 1)$. The subspace $Y_n\cap \{z_1 + \cdots + z_{n} = 0\}$ 
is a deformation retract of the space $Y_n$ and a Stein 
manifold of complex dimension $n-1$. Hence the cohomological 
dimension of the group $P_n$, $\mbox{cd} P_n$, 
is not greater than $n-1$. 
Let $A^*(Y_n)$ be the algebra of all the complex-valued 
differential forms on the space $Y_n$. 
As was shown by Arnol'd \cite{Ar}, the $\bZ$--subalgebra 
generated by the $1$--forms 
$$
\omega_{i, j} := {\frac1{2\pi\sqrt{-1}}}\frac{dz_i - dz_j}{z_i - z_j}, 
\quad 1 \leq i < j \leq n,
$$
is isomorphic to the cohomology algebra $H^*(Y_n; \bZ) = 
H^*(P_n; \bZ)$. Especially in the case $* =1$, $\{[\omega_{i, j}]\}_{1
\leq i<j\leq n}$  is a $\bZ$--free basis of $H^1(P_n; \bZ)$, so that 
$\{[A_{i, j}]\}_{1 \leq
i<j\leq n}$  is a $\bZ$--free basis of $H_1(P_n; \bZ) = {P_n}\abel$.

\begin{lem}\label{2.3}\

\begin{enumerate}
\item ${k_0}^n \neq 0 \in H^n(Y_{n+1}; \La^n\HQ)$, where 
$P_{n+1} = \pi_1(Y_{n+1})$ is regarded as a subgroup of $\bbn =
\pi_1(Y_{n+1}/\sn)$.
\item $h_{n-1}\neq 0 \in H^{n-1}(P_n; \HQ^*\otimes\La^n\HQ)$. 
\end{enumerate}
\end{lem}
\begin{proof}
(1)\qua From \eqref{2-3} and \eqref{2-4} we have 
$$
k_0(A_{i, j}) = \begin{cases}
0, & \mbox{if $i < j \leq n$,}\\
X_i, & \mbox{if $i < j = n+1$,}
\end{cases}
$$
that is 
$$
k_0 = {\sum}^n_{i=1}\omega_{i, n+1}\otimes X_i \in H^1(Y_{n+1}; H).
$$
If we restrict the $n$--form
$$\omega_{1, n+1}\omega_{2, n+1}\cdots\omega_{n, n+1}
= (1/{2\pi\sqrt{-1}})^n{\prod}^n_{i=1}(dz_i -
dz_{n+1}) /(z_i - z_{n+1})$$
to the subspace $Y_{n+1}\cap \{z_{n+1} = 0\}$, then we obtain 
the non-zero $n$--form
$(1/{2\pi\sqrt{-1}})^n$\break ${\prod}^n_{i=1}(dz_i/z_i)$.
Hence the cohomology class 
$$
{k_0}^n = n!\,\omega_{1, n+1}\omega_{2, n+1}\cdots\omega_{n,
n+1}X_1\wedge X_2\wedge \cdots\wedge X_n \in H^n(Y_{n+1}; \La^n\HQ)
$$
does not vanish, as was to be shown.\par
(2)\qua Since $\mbox{cd} P_n \leq n-1$, the Gysin map of the 
extension $$
\Fn \mathop{\to}^{\iota} P_{n+1}\mathop{\to}^{\pi} P_n
$$ 
gives an isomorphism
$$
\pi_\sharp\co  H^n(P_{n+1}; M) \mathop{\to}^{\cong} H^{n-1}(P_n; 
H^*\otimes M)
$$
for any $P_n$--module $M$. Hence $h_{n-1} = \pi_\sharp {k_0}^n 
\neq 0$ by (1).
\end{proof}

The map $r_n\co  \HQ^*\otimes \La^n\HQ \to \La^{n-1}\HQ$ is an isomorphism 
because $\dim_\bQ\HQ = n$. Hence we obtain:
\begin{lem}\label{2.4}
$$
\h_{n-1} \neq 0 \in H^{n-1}(P_n; \La^{n-1}\HQ).
$$
\end{lem}

\section[Proof of \ref{0.1}]{Proof of \fullref{0.1}}\label{sec3}

Our proof of \fullref{0.1} is based on \fullref{1.2} and 
\fullref{2.4}.
For $q \leq n$ we denote by $\pnq$ the set of all
the non-negative partitions 
${\lambda} = (\lambda_1 \geq \lambda_2 \geq \dots \geq \lambda_{n-q} 
\geq 0)$ of $q$ into $n-q$ parts. 
For ${\lambda} 
= (\lambda_1 \geq \lambda_2 \geq \dots \geq \lambda_{n-q} \geq 0) 
\in \pnq$ we introduce a cohomology class $\h_\lambda$ and 
a subgroup $P_\lambda \subset P_n$ by 
\begin{eqnarray*}
\h_\lambda &:=& \h_{\lambda_1}\h_{\lambda_2}\cdots\h_{\lambda_{n-q}} 
\in H^q(B_n; \La^qH_\bQ) \subset H^q(P_n; \La^qH_\bQ),\quad\text{and}\\
P_\lambda &:=& P_{\lambda_1+1}\times P_{\lambda_2+1}\times \cdots\times 
P_{\lambda_{n-q}+1} \subset P_n,
\end{eqnarray*}
respectively. Here $P_{0+1} = P_1$ is the trivial group $\{1\}$. 
Denote by $\iota_\lambda\co  P_\lambda \hookrightarrow P_n$ the
obvious inclusion map and $\varpi_k\co  P_\lambda \to P_{\lambda_k+1}$
the obvious projection. 
\fullref{0.1} follows from: 
\begin{thm}\label{3.1}
The cohomology classes $\{\h_\lambda; \, \lambda \in \pnq\}$ are 
linearly independent in $H^q(P_n; \La^qH_\bQ)$.
\end{thm}
In fact, when $q \leq n/2$, the set of all the non-negative partitions 
of $q$ into $n-q$ parts does not depend on $n$.\par

Endow the partitions $\pnq$ with the lexicographic
order. For example, 
$(q \geq 0 \geq \dots \geq 0)$ is the maximal partition. 
\fullref{3.1} is reduced to the following

\par
\noindent
{\bf Assertions}\qua For any $\lambda$ and $\mu \in \pnq$ we have:

\begin{enumerate}\renewcommand{\theenumi}{\Alph{enumi}}
\item\label{(A)} ${\iota_\lambda}^*\h_\lambda \neq 0 \in H^q(P_\lambda;
\La^qH_\bQ)$
\item\label{(B)}If $\mu \gneqq \lambda$, then ${\iota_\lambda}^*\h_\mu = 0
\in  H^q(P_\lambda; \La^qH_\bQ)$.
\end{enumerate}
\par
In fact, assume we have a nontrivial linear relation 
$$\sum_{\lambda \in \pnq}c_\lambda\h_\lambda = 0 \in 
H^q(P_n; \La^qH_\bQ).$$
Choose the minimum $\lambda$ satisfying $c_\lambda \neq 0$.
Applying ${\iota_\lambda}^*$ to the relation, we obtain 
$c_\lambda{\iota_\lambda}^*\h_\lambda = 0$ from Assertion \ref{(B)}. 
Assertion \ref{(A)} implies $c_\lambda = 0$, 
which contradicts the choice of $\lambda$.

\noindent
\begin{proof}[Proof of Assertion \ref{(A)}]
Let $b_1 \geq b_2 \geq \cdots \geq b_{\lambda_1} 
> b_{\lambda_1+1} = 0$ be the dual partition of $\lambda$. The number 
of $\lambda_k$'s equal to $p$ is $b_{p} - b_{p+1}$. We abbreviate 
$\h_{p, k} := {\varpi_{k}}^*\h_{p}$. 
Since $\operatorname{cd} P_{\lambda_k+1} \leq \lambda_k$, 
we have $\h_{p, k} = 0$ if $p > \lambda_k$, or 
equivalently, $k > b_p$. Moreover we have 
$\h_{\lambda_k, k}\h_{p, k} = 0$ for any $p \geq 1$
since $H^{\lambda_k+p}(P_{\lambda_k+1}; \La^{\lambda_k+p}H_\bQ) = 0$. 
From \fullref{1.2}
we have
$$
{\iota_\lambda}^*\h_p = {\sum}^{n-q}_{k=1}\h_{p, k} 
\in H^p(P_\lambda; \La^pH),
$$
so that
\begin{eqnarray*}
&& {\iota_\lambda}^*\h_\lambda 
= \prod^{n-q}_{k=1}{\iota_\lambda}^*\h_{\lambda_k} 
= \prod^{\lambda_1}_{p=1}({\iota_\lambda}^*\h_p)^{b_p - b_{p+1}} \\
&=& \prod^{\lambda_1}_{p=1}(\h_{p, 1} + \h_{p, 2} 
+ \cdots + \h_{p, n-q})^{b_p - b_{p+1}}\\
&=&\prod^{\lambda_1}_{p=1}(\h_{p, 1} + \h_{p, 2} 
+ \cdots + \h_{p, b_p})^{b_p - b_{p+1}}
= \prod^{\lambda_1}_{p=1}(\h_{p, {b_{p+1}+1}} + \cdots 
+ \h_{p, b_p})^{b_p - b_{p+1}}\\
&=&\prod^{\lambda_1}_{p=1}(b_p - b_{p+1})!\,\h_{p, {b_{p+1}+1}}\cdots 
\h_{p, b_p}\\
&=& \left(\prod^{\lambda_1}_{p=1}(b_p - b_{p+1})!\right)\,
\h_{\lambda_1, 1}\h_{\lambda_2, 2}\cdots\h_{\lambda_{n-q}, n-q}.
\end{eqnarray*}
Here the fifth equal sign comes from the equation 
$\h_{\lambda_k, k}\h_{p, k} = 0$. 
Clearly $r_\lambda := \prod^{\lambda_1}_{p=1}(b_p - b_{p+1})!$ 
is a positive integer. From \fullref{2.4} and the K\"unneth formula 
$\h_{\lambda_1, 1}\h_{\lambda_2, 2}\cdots\h_{\lambda_{n-q}, n-q} \neq 0 
\in H^q(P_\lambda; \La^qH_\bQ)$. This proves Assertion \ref{(A)}.
\end{proof}

\begin{proof}[Proof of Assertion \ref{(B)}]
Suppose $\mu > \lambda$ with respect to the lexicographic order, namely,
$\mu_1 = \lambda_1 \geq \mu_2 = \lambda_2 \geq \dots \geq \mu_h 
= \lambda_h \geq \mu_{h+1} > \lambda_{h+1}$
for some $h$, $0 \leq h < n-q$. 
Let $\nu := (\nu_1 \geq \nu_2 \geq \dots \geq \nu_{h})$ be the (truncated) 
partition of $q' := \lambda_1 + \lambda_2 + \cdots + \lambda_h$ defined by 
$\nu_k := \lambda_k = \mu_k$, $k \leq h$. From Assertion \ref{(A)}  
$$
{\iota_\lambda}^*(\h_{\mu_1}\h_{\mu_2}\cdots\h_{\mu_h}) 
= r_\nu\h_{\mu_1, 1}\h_{\mu_2, 2}\cdots\h_{\mu_h, h} 
\in H^{q'}(P_\lambda; \La^{q'}H).
$$
In fact, from $\mu_h > \lambda_{h+1}$, we have $\h_{\mu_i, j}
= 0$ if $i < j$. Since $\mu_{h+1} \gneqq \lambda_k$ for any $k \geq
h+1$, we have
$$
{\iota_\lambda}^*(\h_{\mu_1}\cdots\h_{\mu_h}\h_{\mu_{h+1}}) 
= r_\nu\h_{\mu_1, 1}\cdots\h_{\mu_h, h}
(\h_{\mu_{h+1},1} + \cdots + \h_{\mu_{h+1},h}) = 0
$$
Hence ${\iota_\lambda}^*(\h_{\mu}) = 0$, as was to be shown.
\end{proof}
This completes the proof of \fullref{3.1} and \fullref{0.1}.

\section{Concluding remarks}\label{sec4}

We conclude this note by giving some remarks on the twisted 
cohomology
of the automorphism group $\An$ and the braid group $B_n$. \par
The IA--automorphism group $\IAn$ is
defined to be the kernel of the  action of the group $\An$ on the
homology group $H = \Fn\abel$. 
We have an extension of groups $IA_n \to \An \to \operatorname{GL}(H)$. 
The map $\tau^\theta_1$ restricted to
$\IAn$ gives an isomorphism of the abelianization of the group $\IAn$
onto the  module $H^*\otimes\La^2H$
$$
\tau_1\co  \IAn\abel \,\mathop{\to}^{\cong}\, H^*\otimes\La^2H
$$ 
(Cohen and Pakianathan \cite{CP}, Farb \cite{F}, Kawazumi
\cite{Ka4}). Here we embed $\La^2H$ into $H^{\otimes 2}$ by $X_i\wedge
X_j \mapsto X_i\otimes X_j - X_j\otimes X_i$ for $1 \leq i, j\leq
n$. \fullref{2.2} implies $\xi^*\co  H^1(\IAn; \bZ) \to H^1(P_n; \bZ)$ is
surjective.  From the result of Arnol'd \cite{Ar} quoted in
\fullref{sec2}, the cohomology algebra $H^*(P_n; \bZ)$ is generated by
the first cohomology classes. Hence we obtain:
\begin{cor}\label{4.1}
The algebra homomorphism
$$
\xi^*\co  H^*(\IAn; \bZ) \to H^*(P_n; \bZ)
$$
induced by the homomorphism $\xi\co  P_n \to \IAn$ 
is surjective.
\end{cor}

It should be remarked that it does {\it not} imply that
the map $\xi^*\co  H^*(\An; M) \to H^*(B_n; M)$ is
surjective for a $\bQ[\operatorname{GL}(H)]$--module $M$. 
In fact, the quotient groups $\An/\IAn = \operatorname{GL}(H)$ 
and $B_n/P_n = \sn$ differ from each other.\par

Fred Cohen \cite[Lemma 7.2, page 261]{CLM} described the action of the
symmetric group $\sn$ on the integral cohomology of the group $P_n$,
$H^*(P_n; \bZ)$. Later Lehrer and Solomon \cite{LS} gave another
explicit description of the $\bQ[\sn]$--module $H^*(P_n; \bQ)$.
Moreover Cohen \cite[Theorem 3.1, page 225]{CLM} computed the twisted
cohomology $H^*(\bn; H^{\otimes m}\otimes \Bbb F)$ for any field $\Bbb
F$ and any $m \geq 0$.  It would be interesting if one could describe
the submodule of $H^*(\bn; M)$ generated by all the possible algebraic
combinations coming from the twisted Morita--Mumford classes $h_p$'s
in an explicit manner. Here we should remark the $\sn$--invariant
inner product $\cdot\co H\otimes H \to \bZ$ defined by $X_i\cdot X_j =
\delta_{i, j}$, $1 \leq i, j \leq n$, gives a $\bn$--isomorphism $H
\cong H^*$.  \par As was stated in Introduction, the algebra
$H^*(\MCG; \La^*\HQ)$ is stably isomorphic to the polynomial algebra
in the twisted Morita--Mumford classes $m_{i, j}$'s.  The intersection
pairing of the surface $\Sigma_{g, 1}$, $H^{\otimes 2} \to \bZ$, gives
an isomorphism $H \cong H^*$ of $\MCG$--modules, so that the cocycle
$\tau^\theta_1$ restricted to $\MCG$ can be regarded as a cocycle
$\tau^\theta_1\co \MCG \to H^{\otimes 3}$.  As was proved by Kawazumi
and Morita in \cite{KM}, for any twisted Morita--Mumford class
$m_{i,j}$ we have an $\MCG$--homomorphism $C\co (H^{\otimes
3})^{\otimes (2i+j-2)}$ $\to \bZ$ obtained from the intersection
pairing such that $C_*[\tau^\theta_1]^{2i+j-2} = m_{i, j}$.  In other
words, the natural map
$$
((\La^*H^1(\mathcal{I}_{g, 1}; \bQ))\otimes M)^{\operatorname{Sp}(H)} \to 
H^*(\MCG; M)
$$
is stably surjective for any finite dimensional 
$\bQ[\operatorname{Sp}(H)]$--module $M$.  Here $\mathcal{I}_{g, 1}$
is the Torelli group, i.e, the kernel of the action of $\MCG$ 
on the homology $H$. \par
Recently Galatius \cite{G} proved 
the rational reduced cohomology $\widetilde{H}^*(\An; \bQ)$ 
vanishes in a stable range. It would be very interesting 
to know whether a similar result holds also for 
twisted coefficients. 
\begin{expect} For a finite dimensional
$\bQ[\operatorname{GL}(H)]$--module $M$, the natural map 
$$
((\La^*H^1(\IAn; \bQ))\otimes M)^{\operatorname{GL}(H)} \to 
H^*(\An; M)
$$
is surjective in some stable range.
\end{expect}
In the case $M$ is the trivial module $\bQ$, this expectation is
exactly the fact that 
$\widetilde{H}^*(\An; \bQ)$ vanishes in some stable range, which 
Galatius \cite{G} proved. A result of Hatcher and Wahl \cite{HW} tells us 
it holds also for $M = (H^*)^{\otimes m}$ for any $m \geq 1$.

\subsubsection*{Acknowledgements}
The author would like to thank Fred Cohen, Hiroaki Terao, 
Hirofumi Yamada and Youichi Shibukawa 
for inspiring discussions. 
He would also like to thank Fred Cohen (once again), 
Benson Farb, Soren Galatius and Nathalie Wahl for 
giving him information about their own 
published/unpublished works.

\bibliographystyle{gtart}
\bibliography{link}

\end{document}